\documentclass[11pt,a4paper]{article}
\pdfoutput=1 
\usepackage[margin=2.4cm]{geometry}
\usepackage{amsmath,amsthm,amssymb,enumerate, bbm ,graphicx,color,caption,upgreek, float, tikz, wasysym, subcaption,booktabs,longtable, appendix,graphics, pdfpages,rotating,mathtools,textcomp, array, blkarray,setspace}
\usepackage{xcolor,colortbl}

\DeclarePairedDelimiter\ceil{\lceil}{\rceil}
\DeclarePairedDelimiter\floor{\lfloor}{\rfloor}
\usepackage[pdftex]{hyperref}

\definecolor{blauw}{RGB}{61,158,255}
\definecolor{donkerblauw}{RGB}{0,0,255}
\definecolor{donkergroen}{RGB}{46,148,0}
\definecolor{donkerrood}{RGB}{204,0,0}
\newenvironment{speciaalenumerate}{
\begin{enumerate}[(i)]
  \setlength{\itemsep}{1pt}
  \setlength{\parskip}{0pt}
  \setlength{\parsep}{0pt}
}{\end{enumerate}}

\makeatletter 
\newcommand\mynobreakpar{\par\nobreak\@afterheading} 
\makeatother

\makeatletter
\let\@fnsymbol\@arabic
\makeatother

\newcommand{\N}{\mathbb{N}}
\newcommand{\Z}{\mathbb{Z}}

\usepackage[english]{babel}

\newtheorem{theorem}{Theorem}[section]

\newtheorem{claim}[theorem]{Claim}
\newtheorem{proposition}[theorem]{Proposition}

\theoremstyle{definition}

\newtheorem*{examp*}{Example}

\newtheorem{remark}{Remark}[section]

\theoremstyle{plain}
\newfloat{Algorithm}{!hbt}{alg}

\newcounter{thm}[section]

\title{New lower bound on the Shannon capacity of \texorpdfstring{$C_7$}{C7} from  circular~graphs} \date{}
\author{Sven Polak\thanks{Korteweg-De Vries Institute for Mathematics, University of Amsterdam. E-mail: \href{mailto:s.c.polak@uva.nl}{s.c.polak@uva.nl}, \href{mailto:a.schrijver@uva.nl}{a.schrijver@uva.nl}. The research leading to these
results has received funding from the European Research Council under the European Union’s Seventh Framework Programme (FP7/2007-2013) / ERC grant agreement \textnumero 339109.}\,\, and Alexander Schrijver\footnotemark[1]}
\begin{document}
\maketitle
\setcounter{footnote}{1}

\noindent \textbf{Abstract.}  
We give an independent set of size~$367$ in the fifth strong product power of~$C_7$, where~$C_7$ is the cycle on~$7$ vertices. This leads to an improved lower bound on the Shannon capacity of~$C_7$: $\Theta(C_7)\geq 367^{1/5} > 3.2578$. The independent set is found by computer, using the fact that the set~$\{t \cdot (1,7,7^2,7^3,7^4) \,\, | \,\, t \in \Z_{382}\} \subseteq \Z_{382}^5$ is independent in the fifth strong product power of the circular graph~$C_{108,382}$. Here the circular graph~$C_{k,n}$ is the graph with vertex set~$\Z_{n}$, the cyclic group of order~$n$, in which two distinct vertices are adjacent if and only if their distance (mod~$n$) is strictly less than~$k$. 

\,$\phantom{0}$



\section{Introduction}
 For any graph~$G=(V,E)$, let~$G^d$ denote the graph with vertex set~$V^d$ and edges between two distinct vertices~$(u_1,\ldots,u_d)$ and~$(v_1,\ldots,v_d)$ if and only if for all~$i \in \{1,\ldots,d\}$ one has either $u_i=v_i$ or~$u_iv_i \in E$. The graph~$G^d$ is known as the~$d$-th \emph{strong product power} of~$G$. The \emph{Shannon capacity} of~$G$ is
\begin{align}
    \Theta(G):= \sup_{d \in \N} \sqrt[d]{\alpha(G^d)},
\end{align}
where~$\alpha(G^d)$ denotes the maximum cardinality of an independent set in~$G^d$, i.e., a set of vertices no two of which are adjacent.
As~$\alpha(G^{d_1+d_2}) \geq \alpha(G^{d_1})\alpha(G^{d_2})$ for any two positive integers~$d_1$ and~$d_2$, by Fekete's lemma~$\cite{fekete}$ it holds that $\Theta(G) = \lim_{d \to \infty} \sqrt[d]{\alpha(G^d)}$. 

The Shannon capacity was introduced by Shannon~$\cite{shannon}$ and is an important and widely studied  parameter in information theory (see e.g.,~\cite{ alon,bohman, haemers,lovasz, zuiddam}). It is the effective size of an alphabet in an information channel represented by the graph~$G$. The input is a set of letters~$V=\{0,\ldots,n-1\}$ and two letters are confusable when transmitted over the channel if and only if there is an edge between them in~$G$. Then~$\alpha(G)$ is the maximum size of a set of pairwise non-confusable single letters. Moreover,~$\alpha(G^d)$ is the maximum size of a set of pairwise non-confusable~$d$-letter words. Taking~$d$-th roots and letting~$d$ go to infinity, we find the effective size of the alphabet in the information channel:~$\Theta(G)$. 

The Shannon capacity of~$C_5$, the cycle on~$5$ vertices, was already discussed by Shannon in 1956~$\cite{shannon}$. It was determined more than twenty years later by Lov\'asz~$\cite{lovasz}$ using his famous~$\vartheta$-function. He proved that~$\Theta(C_5)= \sqrt{5}$. The  easy lower bound is obtained from the independent set $\{(0,0),(1,2),(2,4),(3,1),(4,3)\}$ in~$C_5^2$ and the ingenious upper bound is given by Lov\'asz's $\vartheta$-function.  More generally, for odd~$n$,
\begin{align} \label{vartheta}
    \Theta(C_n) \leq \vartheta(C_n) = \frac{n\cos(\pi/n)}{1+\cos(\pi/n)}.
\end{align}
For~$n$ even it is not hard to see that~$\Theta(C_n)=n/2$. 

The Shannon capacity of~$C_7$ is still unknown and its determination is a notorious open problem in extremal combinatorics~\cite{bohman, godsil}. Many lower bounds have been given by explicit independent sets in some fixed power of~$C_7$ \cite{baumert, matos, veszer}, while the best known upper bound is~$\Theta(C_7)\leq \vartheta(C_7) < 3.3177$. Here we give an independent set of size~$367$ in~$C_7^5$, which yields~$\Theta(C_7)\geq 367^{1/5} > 3.2578$. The best previously known lower bound on~$\Theta(C_7)$ is~$\Theta(C_7) \geq 350^{1/5} > 3.2271$, found by Mathew and \"Osterg{\aa}rd~$\cite{matos}$. They proved that~$\alpha(C_7^5) \geq 350$ using stochastic search methods that utilize the symmetry of the problem. In~$\cite{baumert}$, a construction is given of an independent set of size~$7^3=343$ in~$C_7^5$. The best known lower bound on~$\alpha(C_7^4)$ is~$108$, by Vesel and \v{Z}erovnik~\cite{veszer}. See Table~\ref{knownindep} for the currently best known bounds on~$\alpha(C_7^d)$ for small~$d$.

\begin{table}[ht]
\centering
\begin{tabular}[t]{l|ccccc}
$d$ & 1 & 2 & 3 & 4 & 5 \\ \hline 
$\alpha(C_7^d)$ & 3 & $10^a$ & $33^d$ & $108^e$--$115^b$ & $367^f$--$401^c$
\end{tabular} 
\caption{Bounds on~$\alpha(C_7^d)$. {\small Key:
\\ $^a$ $\alpha(C_n^2)= \floor{(n^2-n)/4}$ \cite[Theorem 2]{baumert} 
\\$^b$ $\alpha(C_n^d) \leq \alpha(C_n^{d-1})n/2$ \cite[Lemma 2]{baumert} 
\\$^c$ $\alpha(G^d) \leq \vartheta(G)^d$ by Lov\'asz $\cite{lovasz}$
\\$^d$ Baumert et al.~\cite{baumert}
\\$^e$ Vesel and \v{Z}erovnik~\cite{veszer}
\\$^f$ this paper, see the Appendix for the explicit independent set.}} \label{knownindep}
\end{table}
For comparison,~$\alpha(C_7^3)^{1/3} =33^{1/3} \approx 3.2075$, $\alpha(C_7^4)^{1/4} \geq 108^{1/4} \approx 3.2237$ and the previously best known lower bound on~$\alpha(C_7^5)^{1/5}$  is~$350^{1/5} \approx 3.2271$. Now we know that~$\alpha(C_7^5)\geq 367 > 3.2578^5$.

The paper is organized as follows. In Section~\ref{circular} we will examine the circular graphs~$C_{k,n}$. We give a construction that yields independent sets in certain~$C_{k,n}^d$, and we give an explicit description of an independent set~$S$ of size~$382$ in the graph~$C_{108,382}^5$. This independent set does not translate directly to an independent set in~$C_7^5$. However, in Section~$\ref{descr}$ we describe how one can obtain an independent set of size~$367$ in~$C_7^5$ from~$S$, by adapting~$S$, removing vertices and adding new ones. This independent set  is given explicitly in the Appendix.

\section{Circular Graphs \label{circular}}

For  two integers~$a,b$, let~$[a,b]$ denote the set~$\{a,a+1,\ldots,b\}$. For~$k,n \in \Z$ with~$n \geq 2k$, the \emph{circular graph} $C_{k,n}$ is the graph with vertex set~$\Z_n$, the cyclic group of order~$n$, in which two distinct vertices are adjacent if and only if their distance (mod~$n$) is strictly less than~$k$. In other words, it is the \emph{circulant graph} on~$\Z_n$ with \emph{generating set}~$[1,k-1]$, which means that~$V=\Z_n$ and~$E=\{ \{u,v\}\,\, | \,\, u-v \in [1,k-1]\}$. So~$C_{2,n} = C_n$ and~$C_{k,n}$ has vertex set~$V(C_n)=\Z_n$. A closed formula for~$\vartheta(C_{k,n})$, Lov\'asz's upper bound on $\Theta(C_{k,n})$, is given in~$\cite{bachoc}$.

Note that by definition there is an edge between two distinct vertices~$x,y$ of~$C_{k,n}^d$ if and only if there is an~$i\in [1,d] $ such that~$x_i-y_i \pmod{n}$ is either strictly smaller than~$k$ or strictly larger than~$n-k$.  For distinct~$u,v$ in~$\Z_n^d$, define their \emph{distance} to be the maximum over the distances of~$u_i$ and~$v_i$ (mod~$n$), where~$i$ ranges from~$1$ to~$d$. The \emph{minimum distance}~$d_{\text{min}}(D)$ of a set~$D \subseteq \Z_n^d$ is the minimum distance  between any pair of distinct elements of~$D$. (If~$|D|=1$, set~$d_{\text{min}}(D)=\infty$.) Then~$d_{\text{min}}(D) \geq k$ if and only if~$D$ is independent in~$C_{k,n}^d$.

A homomorphism from a graph $G_1 = (V_1,E_1)$ to a graph $G_2 = (V_2,E_2)$ is a function $f\, :\, V_1 \to V_2$ such that if~$ij \in E_1$
then~$f(i)f(j) \in E_2$ (in particular,~$f(i) \neq f(j)$). If there exists a homomorphism~$f \, : \, G_1 \to G_2$ we write~$G_1 \to G_2$. For any graph~$G$, we write~$\overline{G}$ for the complement of~$G$. If~$\overline{G} \to \overline{H}$, then~$\alpha(G) \leq \alpha(H)$ and~$\Theta(G) \leq \Theta(H)$.  The circular graphs have the property that~$\overline{C_{k',n'}} \to \overline{C_{k,n}}$ if and only if~$n'/k' \leq n/k$~\cite{starnote}. So if~$n'/k' \leq n/k$, then~$\alpha(C_{k',n'}^d) \leq \alpha(C_{k,n}^d)$ (for any~$d$) and~$\Theta(C_{k',n'}) \leq \Theta(C_{k,n})$. Moreover,~$\alpha(C_{k,n}^d)$ and~$\Theta(C_{k,n})$ only depend on the fraction~$n/k$.

An independent set in~$C_{k,n}^d$ gives an independent set in~$C_{\ceil{2n/k}}^d$, since $\overline{C_{k,n}} \to \overline{C_{2,\ceil{2n/k}}}$. Explicitly, consider the elements of~$\Z_n$ as integers between~$0$ and~$n-1$ and replace each element~$i$ by~$\floor{2i/k}$, and consider the outcome as an element of~$\Z_{\ceil{2n/k}}$. This gives indeed a homomorphism $\overline{C_{k,n}} \to \overline{C_{2,\ceil{2n/k}}}$ as the image of any two elements with distance at least~$k$ has distance at least~$2$.
 
 First, we will give an independent set of size~$382$ in~$C_{108,382}^5$. As~$382/108 > 7/2$ this does not directly give an independent set in~$C_{2,7}^5$. However, in Section~\ref{descr} we try to adapt the independent set, remove some words and add as many new words as possible to obtain a large independent set in~$C_7^5$. 

\begin{proposition}\label{is}
The set~$S:=\{t \cdot (1,7,7^2,7^3,7^4) \,\, | \,\, t \in \Z_{382}\} \subseteq \Z_{382}^5$ is independent in $C_{108,382}^5$. 
\end{proposition}
\proof
If~$x,y \in S$ then also~$x-y \in S$. So it suffices to check that for all nonzero~$x \in S$:
\begin{align} \label{coor108}
\text{$\exists \,\, i \in [1,5]$ such that~$x_i \in [108, 274]$.}
\end{align}
Let~$x=t \cdot (1,7,7^2,7^3,7^4) \in S$ be arbitrary, with~$0 \neq t \in \Z_{382}$.  For~$t \in [108,274]$ clearly~$\eqref{coor108}$ holds with~$i=1$ (as then~$x_i =t \in [108,274]$).  Also we have  $[275,381]=-[1,107]$, so it suffices to verify~\eqref{coor108} for~$t \in [1,107]$. Note that for~$t \in [16,39]$ one has~$108\leq 7t\leq274 $, so~\eqref{coor108} is satisfied with~$i=2$. Also note that~$69\cdot 7 \equiv 101 \pmod{382}$. So for~$t \in [70,93]$ one has~$7t \equiv 101+7(t-69) \pmod{382} \in [108, 274]$, i.e.,~\eqref{coor108} is satisfied with~$i=2$. For the remaining~$t \in [1,107]$, please take a glance at Table~\ref{verif}. In each row, in each of the three subtables, there is at least one entry in~$[108,274]$. This completes the proof. 
\endproof 

\begin{table}[ht]
\centering
\tiny
\begin{tabular}[t]{rrrrr}
1 & 7 & 49 & 343 & \cellcolor{donkergroen!25} 109 \\ 
2 & 14 & 98 & 304 & \cellcolor{donkergroen!25} 218 \\ 
3 & 21 & \cellcolor{donkergroen!25}147 & \cellcolor{donkergroen!25} 265 & 327 \\ 
4 & 28 & \cellcolor{donkergroen!25}196 & \cellcolor{donkergroen!25}226 & 54 \\ 
5 & 35 & \cellcolor{donkergroen!25}245 & \cellcolor{donkergroen!25}187 &\cellcolor{donkergroen!25} 163 \\ 
6 & 42 & 294 &\cellcolor{donkergroen!25} 148 &\cellcolor{donkergroen!25} 272 \\ 
7 & 49 & 343 &\cellcolor{donkergroen!25} 109 & 381 \\ 
8 & 56 & 10 & 70 &\cellcolor{donkergroen!25} 108 \\ 
9 & 63 & 59 & 31 &\cellcolor{donkergroen!25} 217 \\ 
10 & 70 &\cellcolor{donkergroen!25} 108 & 374 & 326 \\ 
11 & 77 &\cellcolor{donkergroen!25} 157 & 335 & 53 \\ 
12 & 84 &\cellcolor{donkergroen!25} 206 & 296 &\cellcolor{donkergroen!25} 162 \\ 
13 & 91 &\cellcolor{donkergroen!25} 255 &\cellcolor{donkergroen!25} 257 &\cellcolor{donkergroen!25} 271 \\ 
14 & 98 & 304 &\cellcolor{donkergroen!25} 218 & 380 \\ 
15 & 105 & 353 &\cellcolor{donkergroen!25} 179 & 107 \\ 
&&&& \\ 
40 & 280 & 50 & 350 &\cellcolor{donkergroen!25} 158 \\ 
41 & 287 & 99 & 311 &\cellcolor{donkergroen!25} 267 \\ 
42 & 294 &\cellcolor{donkergroen!25} 148 &\cellcolor{donkergroen!25} 272 & 376 \\ 
43 & 301 & \cellcolor{donkergroen!25}197 & \cellcolor{donkergroen!25}233 & 103 \\ 
44 & 308 &\cellcolor{donkergroen!25} 246 & \cellcolor{donkergroen!25}194 & \cellcolor{donkergroen!25}212 \\ 
\end{tabular}\,\,\,\,\,\,\,\,\,\,\,\,\,\,\,
\begin{tabular}[t]{rrrrr}
45 & 315 & 295 &\cellcolor{donkergroen!25} 155 & 321 \\ 
46 & 322 & 344 & \cellcolor{donkergroen!25}116 & 48 \\ 
47 & 329 & 11 & 77 & \cellcolor{donkergroen!25}157 \\ 
48 & 336 & 60 & 38 &\cellcolor{donkergroen!25} 266 \\ 
49 & 343 &\cellcolor{donkergroen!25} 109 & 381 & 375 \\ 
50 & 350 & \cellcolor{donkergroen!25}158 & 342 & 102 \\ 
51 & 357 & \cellcolor{donkergroen!25}207 & 303 &\cellcolor{donkergroen!25} 211 \\ 
52 & 364 &\cellcolor{donkergroen!25} 256 & \cellcolor{donkergroen!25}264 & 320 \\ 
53 & 371 & 305 &\cellcolor{donkergroen!25} 225 & 47 \\ 
54 & 378 & 354 &\cellcolor{donkergroen!25} 186 &\cellcolor{donkergroen!25} 156 \\ 
55 & 3 & 21 &\cellcolor{donkergroen!25} 147 &\cellcolor{donkergroen!25} 265 \\ 
56 & 10 & 70 & \cellcolor{donkergroen!25}108 & 374 \\ 
57 & 17 &\cellcolor{donkergroen!25} 119 & 69 & 101 \\ 
58 & 24 & \cellcolor{donkergroen!25}168 & 30 &\cellcolor{donkergroen!25} 210 \\ 
59 & 31 &\cellcolor{donkergroen!25} 217 & 373 & 319 \\ 
60 & 38 &\cellcolor{donkergroen!25} 266 & 334 & 46 \\ 
61 & 45 & 315 & 295 &\cellcolor{donkergroen!25} 155 \\ 
62 & 52 & 364 & \cellcolor{donkergroen!25}256 &\cellcolor{donkergroen!25} 264 \\ 
63 & 59 & 31 & \cellcolor{donkergroen!25}217 & 373 \\ 
64 & 66 & 80 &\cellcolor{donkergroen!25} 178 & 100 \\ 
\end{tabular} \,\,\,\,\,\,\,\,\,\,\,\,\,\,\,
\begin{tabular}[t]{rrrrr}
65 & 73 &\cellcolor{donkergroen!25} 129 &\cellcolor{donkergroen!25} 139 &\cellcolor{donkergroen!25} 209 \\ 
66 & 80 &\cellcolor{donkergroen!25} 178 & 100 & 318 \\ 
67 & 87 &\cellcolor{donkergroen!25} 227 & 61 & 45 \\ 
68 & 94 & 276 & 22 &\cellcolor{donkergroen!25} 154 \\ 
69 & 101 & 325 & 365 &\cellcolor{donkergroen!25} 263 \\ 
&&&& \\ 
94 & 276 & 22 & \cellcolor{donkergroen!25}154 & 314 \\ 
95 & 283 & 71 & \cellcolor{donkergroen!25}115 & 41 \\ 
96 & 290 &\cellcolor{donkergroen!25} 120 & 76 & \cellcolor{donkergroen!25}150 \\ 
97 & 297 & \cellcolor{donkergroen!25}169 & 37 &\cellcolor{donkergroen!25} 259 \\ 
98 & 304 &\cellcolor{donkergroen!25} 218 & 380 & 368 \\ 
99 & 311 & \cellcolor{donkergroen!25}267 & 341 & 95 \\ 
100 & 318 & 316 & 302 & \cellcolor{donkergroen!25}204 \\ 
101 & 325 & 365 &\cellcolor{donkergroen!25} 263 & 313 \\ 
102 & 332 & 32 &\cellcolor{donkergroen!25} 224 & 40 \\ 
103 & 339 & 81 & \cellcolor{donkergroen!25}185 & \cellcolor{donkergroen!25}149 \\ 
104 & 346 & \cellcolor{donkergroen!25}130 & \cellcolor{donkergroen!25}146 &\cellcolor{donkergroen!25} 258 \\ 
105 & 353 & \cellcolor{donkergroen!25}179 & 107 & 367 \\ 
106 & 360 & \cellcolor{donkergroen!25}228 & 68 & 94 \\ 
107 & 367 & 277 & 29 & \cellcolor{donkergroen!25}203 \\ 
\end{tabular}
\caption{Part of the verification that~$S$ is independent in~$C_{108,382}^5$.} \label{verif}
\end{table}
\noindent  The authors found the above independent set when looking for answers to the following question.
\begin{align} \label{question}
  \text{For~$n,d,q$, what is the minimum distance~$k(n,d,q)$ of }   \{ t \cdot (1,q,\ldots, q^{d-1}) \,\, | \,\, t \in \Z_{n} \} \subseteq \Z_n^d \, \text{?}
 \end{align} 
 The independent set from Proposition~\ref{is} was found by computer (with~$n\geq 350$ and~$d=5$ such that~$n/k(n,d,q)$ is close to~$7/2$). Question~$(\ref{question})$ seems not easy to solve in general.

\section{Description of the method \label{descr}}
Here we describe how to use the independent set from Proposition~\ref{is} to find an independent set of size~$367$ in~$C_7^5$. The procedure is as follows.
\begin{speciaalenumerate}
    \item Start with the independent set~$S$ in~$C_{108,382}^5$ from Proposition~\ref{is}.
    \item Add the word~$(40,123,40,123,40)$ mod~$382$ to each word in~$S$.
    \item Replace each letter~$i$, which we now consider to be an integer between~$0$ and~$381$ and not anymore an element in~$\Z_{382}$, in each word from~$S$ by~$\floor{i/54.5}$. Now we have a set of words~$S'$ with only symbols in~$[0,6]$ in it, which we consider as elements of~$\Z_7$.
    \item Remove each word~$u\in S'$ for which there is a~$v \in S'$ such that~$uv \in E(C_7^5)$ from~$S'$, i.e., we remove~$u$ if there is a~$v \in S'$ with~$v \neq u$ such that~$u_i-v_i\in \{0,1,6\}$ for all~$i \in [1,5]$. We denote the set of words which are not removed from~$S'$ by this procedure by~$M$. The computer finds~$|M| = 327$. Note that~$M$ is independent in~$C_7^5$.
    \item \label{v} Find the best possible extension of~$M$ to a larger independent set in~$C_7^5$. To do this, consider the subgraph~$G$ of~$C_7^5$ induced by the words~$x$ in~$\Z_7^5$ with the property that~$M\cup \{x\}$ is independent in~$C_7^5$. This graph is not large, in this case it has~$71$ vertices and~$85$ edges, so a computer finds a maximum size independent set~$I$ in~$G$ quickly. The computer finds $|I|=\alpha(G)=40$, so we can add~$40$ words to~$M$. Write~$R:=M \cup I$. Then~$|R|=327+40=367$ and~$R$ is independent in~$C_7^5$.
\end{speciaalenumerate}
The maximum size independent set~$I$ in the graph~$G$ in (\ref{v}) was found using Gurobi~$\cite{gurobi}$. In steps (ii) and (iii), many possibilities for adding a constant word and for the division factor were tried, but no independent set of size~$368$ or larger was found. Also, the independent set~$R$ of size~$367$ did not seem to be easily extendable. A local search was performed, showing that there exists no triple of words from~$R$ such that if one removes these three words from~$R$, four words can be added to obtain an independent set of size~$368$ in~$C_7^5$. 

 \begin{remark}
One other new bound on~$\alpha(C_{n}^d)$  was obtained (for~$n \leq 15$ and~$d \leq 5$) using independent sets of the form from~\eqref{question}. With $n= 4009$, $d=5$ and $q=27$, we found~$k(n,d,q)= 729$. As~$n/(k(n,d,q))=4009/729<11/2$, this directly yields the new lower bound~$\alpha(C_{11}^5) \geq 4009$. The previously best known lower bound is~$\alpha(C_{11}^5)\geq 3996$ from~\cite{circular, matos}. However, the new lower bound on~$\alpha(C_{11}^5)$ does not imply a new lower bound on~$\Theta(C_{11})$. It is known that $\Theta(C_{11}) \geq \alpha(C_{11}^3)^{1/3} = 148^{1/3}>5.2895$ (cf.~\cite{baumert}), which is larger than~$4009^{1/5}$. 
\end{remark}

\section*{Appendix: explicit code}
The following~$367$ words form an independent set in~$C_7^5$, which proves the new bound~$\Theta(C_7) \geq 367^{1/5} > 3.2578$. It is the set~$R$ from Section~$\ref{descr}$. 

$\phantom{,}$
\begin{samepage}
{\begin{spacing}{0.1}\spaceskip=0.0875em\tiny\noindent \texttt{02020, 02112, 02204, 02306, 02461, 02553, 03645, 03040, 03032, 03124, 03226, 03311, 03403, 14144, 14231, 14323, 14415, 14510, 15602, 15064, 15166, 15251, 15343, 15430, 15522, 16614, 16016, 16101, 16263, 16355, 16450, 16542, 10636, 10021, 10113, 10205, 10300, 10462, 10554, 11656, 11041, 11033, 11125, 11220, 11312, 11404, 11506, 12661, 12053, 12145, 12240, 12232, 12324, 12426, 12511, 13603, 13065, 13160, 13252, 13344, 13446, 13431, 24010, 24102, 24264, 24366, 24451, 24543, 25630, 25022, 25114, 25216, 25301, 25463, 25555, 26650, 26042, 26034, 26136, 26221, 26313, 26405, 26500, 20662, 20054, 20156, 20241, 20233, 20325, 20420, 20512, 21604, 21006, 21161, 21253, 21345, 21440, 21432, 22626, 22011, 22103, 22265, 22360, 22452, 22544, 23631, 23023, 23115, 23210, 23302, 23464, 23566, 34130, 34222, 34314, 34416, 34501, 35663, 35055, 35150, 35242, 35234, 35336, 35421, 35513, 36605, 36000, 36162, 36254, 36356, 36441, 36433, 30620, 30012, 30104, 30206, 30361, 30453, 30545, 31632, 31024, 31126, 31211, 31303, 31465, 31560, 32652, 32044, 32131, 32223, 32315, 32410, 32502, 33664, 33066, 33151, 33243, 33235, 33330, 33422, 44616, 44001, 44163, 44255, 44350, 44442, 44434, 44536, 45621, 45013, 45105, 45200, 45362, 45454, 45556, 46633, 46025, 46120, 46212, 46304, 46406, 46561, 40653, 40045, 40132, 40224, 40326, 40411, 40503, 41665, 41060, 41152, 41244, 41331, 41423, 41515, 42610, 42002, 42164, 42266, 42351, 42443, 42435, 43622, 43014, 43116, 43201, 43363, 43455, 43550, 54634, 54036, 54121, 54213, 54305, 54400, 54562, 55654, 55056, 55141, 55133, 55225, 55320, 55412, 55504, 56606, 56061, 56153, 56245, 56332, 56424, 56526, 50611, 50003, 50165, 50260, 50352, 50444, 51623, 51015, 51110, 51202, 51364, 51551, 52643, 52635, 52030, 52122, 52214, 53655, 53134, 64332, 64424, 64526, 65611, 65003, 65260, 65352, 65444, 65546, 66623, 66110, 66202, 66364, 66466, 66551, 60643, 60645, 60030, 60122, 60214, 60316, 60401, 60563, 61050, 61142, 61134, 61236, 61321, 61413, 62600, 62062, 62154, 62256, 62341, 62333, 62520, 63612, 63004, 63106, 63261, 63353, 63445, 63540, 64532, 04026, 04111, 04203, 04460, 04552, 05644, 05031, 05123, 05310, 05402, 05564, 06666, 06051, 06143, 06230, 06322, 06414, 06516, 00601, 00063, 00155, 00250, 00342, 00334, 00436, 01613, 01100, 01262, 01354, 01456, 01541, 02625, 00521, 01005, 02533, 03565, 04052, 04365, 04624, 04660, 05046, 05225, 10534, 14246, 15435, 22524, 24615, 24651, 32046, 34035, 34043, 36525, 40040, 41246, 42530, 43514, 45641, 50531, 51456, 52400, 52563, 53050, 53142, 53320, 53412, 56340, 61505, 62425, 64154, 64340, 65105, 66025}. \end{spacing}}%
\end{samepage}

\section*{Acknowledgements}
The authors want to thank Bart Litjens, Bart Sevenster, Jeroen Zuiddam and the two anonymous referees for very useful comments.

\selectlanguage{english}

\end{document}